# REJOINDER TO "LEAST ANGLE REGRESSION" BY EFRON ET AL.

By Bradley Efron, Trevor Hastie, Iain Johnstone and Robert Tibshirani

The original goal of this project was to explain the striking similarities between models produced by the Lasso and Forward Stagewise algorithms, as exemplified by Figure 1. LARS, the Least Angle Regression algorithm, provided the explanation and proved attractive in its own right, its simple structure permitting theoretical insight into all three methods. In what follows "LAR" will refer to the basic, unmodified form of Least Angle Regression developed in Section 2, while "LARS" is the more general version giving LAR, Lasso, Forward Stagewise and other variants as in Section 3.4. Here is a summary of the principal properties developed in the paper:

1. LAR builds a regression model in piecewise linear forward steps, accruing explanatory variables one at a time; each step is taken along the equiangular direction between the current set of explanators. The step size is less greedy than classical forward stepwise regression, smoothly blending in new variables rather than adding them discontinuously.
2. Simple modifications of the LAR procedure produce all Lasso and Forward Stagewise solutions, allowing their efficient computation and showing that these methods also follow piecewise linear equiangular paths. The Forward Stagewise connection suggests that LARS-type methods may also be useful in more general "boosting" applications.
3. The LARS algorithm is computationally efficient; calculating the full set of LARS models requires the same order of computation as ordinary least squares.
4. A $k$-step LAR fit uses approximately $k$ degrees of freedom, in the sense of added prediction error (4.5). This approximation is exact in the case of orthogonal predictors and is generally quite accurate. It permits $C_p$-type stopping rules that do not require auxiliary bootstrap or cross-validation computations.
5. For orthogonal designs, LARS models amount to a succession of soft thresholding estimates, (4.17).

All of this is rather technical in nature, showing how one might efficiently carry out a program of automatic model-building ("machine learning"). Such programs seem increasingly necessary in a scientific world awash in huge data sets having hundreds or even thousands of available explanatory variables.

What this paper, strikingly, does not do is justify any of the three algorithms as providing *good* estimators in some decision-theoretic sense. A few hints appear, as in the simulation study of Section 3.3, but mainly we are relying on recent literature to say that LARS methods are at least reasonable algorithms and that it is worthwhile understanding their properties. Model selection, the great underdeveloped region of classical statistics, deserves careful theoretical examination but that does not happen here. We are not as pessimistic as Sandy Weisberg about the potential of automatic model selection, but agree that it requires critical examination as well as (over) enthusiastic algorithm building.

The LARS algorithm in any of its forms produces a one-dimensional path of prediction vectors going from the origin to the full least-squares solution. (Figures 1 and 3 display the paths for the diabetes data.) In the LAR case we can label the predictors $\widehat{\boldsymbol{\mu}}(k)$, where $k$ is identified with both the number of steps and the degrees of freedom. What the figures do not show is when to stop the model-building process and report $\widehat{\boldsymbol{\mu}}$ back to the investigator. The examples in our paper rather casually used stopping rules based on minimization of the $C_p$ error prediction formula.

Robert Stine and Hemant Ishwaran raise some reasonable doubts about $C_p$ minimization as an effective stopping rule. For any one value of $k$, $C_p$ is an unbiased estimator of prediction error, so in a crude sense $C_p$ minimization is trying to be an unbiased estimator of the optimal stopping point $k_{\text{opt}}$. As such it is bound to overestimate $k_{\text{opt}}$ in a large percentage of the cases, perhaps near 100% if $k_{\text{opt}}$ is near zero.

We can try to improve $C_p$ by increasing the *df* multiplier "2" in (4.5). Suppose we change 2 to some value *mult*. In standard normal-theory model building situations, for instance choosing between linear, quadratic, cubic, ... regression models, the *mult* rule will prefer model $k+1$ to model $k$ if the relevant $t$-statistic exceeds $\sqrt{mult}$ in absolute value (here we are assuming $\sigma^2$ known); $mult = 2$ amounts to using a rejection rule with $\alpha = 16\%$. Stine's interesting $S_p$ method chooses *mult* closer to 4, $\alpha = 5\%$.

This works fine for Stine's examples, where $k_{\text{opt}}$ is indeed close to zero. We tried it on the simulation example of Section 3.3. Increasing *mult* from 2 to 4 decreased the average selected step size from 31 to 15.5, but with a small increase in actual squared estimation error. Perhaps this can be taken as support for Ishwaran's point that since LARS estimates have a broad plateau of good behavior, one can often get by with much smaller models than suggested by $C_p$ minimization. Of course no one example is conclusive in an area as multifaceted as model selection, and perhaps no 50 examples either. A more powerful theory of model selection is sorely needed, but until it comes along we will have to make do with simulations, examples and bits and pieces of theory of the type presented here.

Bayesian analysis of prediction problems tends to favor *much* bigger choices of *mult*. In particular the Bayesian information criterion (BIC) uses $mult = \log(\text{sample size})$. This choice has favorable consistency properties, selecting the correct model with probability 1 as the sample size goes to infinity. However, it can easily select too-small models in nonasymptotic situations.

Jean-Michel Loubes and Pascal Massart provide two interpretations using penalized estimation criteria in the orthogonal regression setting. The first uses the link between soft thresholding and $\ell_1$ penalties to motivate entropy methods for asymptotic analysis. The second is a striking perspective on the use of $C_p$ with LARS. Their analysis suggests that our usual intuition about $C_p$, derived from selecting among projection estimates of different ranks, may be misleading in studying a nonlinear method like LARS that combines thresholding and shrinkage. They rewrite the LARS-$C_p$ expression (4.5) in terms of a penalized criterion for selecting among orthogonal projections. Viewed in this unusual way (for the estimator to be used is *not* a projection!), they argue that *mult* in fact behaves like $\log(n/k)$ rather than 2 (in the case of a $k$-dimensional projection). It is indeed remarkable that this same model-dependent value of *mult*, which has emerged in several recent studies [Foster and Stine (1997), George and Foster (2000), Abramovich, Benjamini, Donoho and Johnstone (2000) and Birgé and Massart (2001)], should also appear as relevant for the analysis of LARS. We look forward to the further extension of the Birgé–Massart approach to handling these nondeterministic penalties.

Cross-validation is a nearly unbiased estimator of prediction error and as such will perform similarly to $C_p$ (with $mult = 2$). The differences between the two methods concern generality, efficiency and computational ease. Cross-validation, and nonparametric bootstrap methods such as the 632+ rule, can be applied to almost any prediction problem. $C_p$ is more specialized, but when it does apply it gives more efficient estimates of prediction error [Efron (2004)] at almost no computational cost. It applies here to LAR, at least when $m < n$, as in David Madigan and Greg Ridgeway's example.

We agree with Madigan and Ridgeway that our new LARS algorithm may provide a boost for the Lasso, making it more useful and attractive for data analysts. Their suggested extension of LARS to generalized linear models is interesting. In logistic regression, the $L_1$-constrained solution is not piecewise linear and hence the pathwise optimization is more difficult. Madigan and Ridgeway also compare LAR and Lasso to least squares boosting for prediction accuracy on three real examples, with no one method prevailing.

Saharon Rosset and Ji Zhu characterize a class of problems for which the coefficient paths, like those in this paper, are piecewise linear. This is a useful advance, as demonstrated with their robust version of the Lasso, and the $\ell_1$-regularized Support Vector Machine. The former addresses some of the robustness concerns of Weisberg. They also report on their work that strengthens the connections between $\varepsilon$-boosting and $\ell_1$-regularized function fitting.

Berwin Turlach's example with uniform predictors surprised us as well. It turns out that 10-fold cross-validation selects the model with $|\beta_1| \approx 45$ in his Figure 3 (left panel), and by then the correct variables are active and the interactions have died down. However, the same problem with 10 times the noise variance does not recover in a similar way. For this example, if the $X_j$ are uniform on $[-\frac{1}{2}, \frac{1}{2}]$ rather than $[0, 1]$, the problem goes away, strongly suggesting that proper centering of predictors (in this case the interactions, since the original variables are automatically centered by the algorithm) is important for LARS.

Turlach also suggests an interesting proposal for enforcing marginality, the hierarchical relationship between the main effects and interactions. In his notation, marginality says that $\beta_{i:j}$ can be nonzero only if $\beta_i$ and $\beta_j$ are nonzero. An alternative approach, more in the "continuous spirit" of the Lasso, would be to include constraints

$$|\beta_{i:j}| \leq \min\{|\beta_i|, |\beta_j|\}.$$

This implies marginality but is stronger. These constraints are linear and, according to Rosset and Zhu above, a LARS-type algorithm should be available for its estimation. Leblanc and Tibshirani (1998) used constraints like these for shrinking classification and regression trees.

As Turlach suggests, there are various ways to restate the LAR algorithm, including the following nonalgebraic purely statistical statement in terms of repeated fitting of the residual vector $\mathbf{r}$:

1. Start with $\mathbf{r} = \mathbf{y}$ and $\widehat{\beta}_j = 0 \ \forall j$.
2. Find the predictor $\mathbf{x}_j$ most correlated with $\mathbf{r}$.
3. Increase $\widehat{\beta}_j$ in the direction of the sign of $\text{corr}(\mathbf{r}, \mathbf{x}_j)$ until some other competitor $\mathbf{x}_k$ has as much correlation with the current residual as does $\mathbf{x}_j$.
4. Update $\mathbf{r}$, and move $(\widehat{\beta}_j, \widehat{\beta}_k)$ in the joint least squares direction for the regression of $\mathbf{r}$ on $(\mathbf{x}_j, \mathbf{x}_k)$ until some other competitor $\mathbf{x}_\ell$ has as much correlation with the current residual.
5. Continue in this way until all predictors have been entered. Stop when $\text{corr}(\mathbf{r}, \mathbf{x}_j) = 0 \ \forall j$, that is, the OLS solution.

Traditional forward stagewise would have completed the least-squares step at each stage; here it would go only a fraction of the way, until the next competitor joins in.

Keith Knight asks whether Forward Stagewise and LAR have implicit criteria that they are optimizing. In unpublished work with Trevor Hastie, Jonathan Taylor and Guenther Walther, we have made progress on that question. It can be shown that the Forward Stagewise procedure does a sequential minimization of the residual sum of squares, subject to

$$\sum_j \left| \int_0^t \beta_j'(s)\, ds \right| \leq t.$$

This quantity is the total $L_1$ arc-length of the coefficient curve $\beta(t)$. If each component $\beta_j(t)$ is monotone nondecreasing or nonincreasing, then $L_1$ arc-length equals the $L_1$-norm $\sum_j |\beta_j|$. Otherwise, they are different and $L_1$ arc-length discourages sign changes in the derivative. That is why the Forward Stagewise solutions tend to have long flat plateaus. We are less sure of the criterion for LAR, but currently believe that it uses a constraint of the form $\sum_j |\int_0^k \beta_j(s)\, ds| \leq A$.

Sandy Weisberg, as a ranking expert on the careful analysis of regression problems, has legitimate grounds for distrusting automatic methods. Only foolhardy statisticians dare to ignore a problem's context. (For instance it helps to know that diabetes progression behaves differently after menopause, implying strong age–sex interactions.) Nevertheless even for a "small" problem like the diabetes investigation there is a limit to how much context the investigator can provide. After that one is drawn to the use of automatic methods, even if the "automatic" part is not encapsulated in a single computer package.

In actual practice, or at least in good actual practice, there is a cycle of activity between the investigator, the statistician and the computer. For a multivariable prediction problem like the diabetes example, LARS-type programs are a good first step toward a solution, but hopefully not the last step. The statistician examines the output critically, as did several of our commentators, discussing the results with the investigator, who may at this point suggest adding or removing explanatory variables, and so on, and so on.

Fully automatic regression algorithms have one notable advantage: they permit an honest evaluation of estimation error. For instance the $C_p$-selected LAR quadratic model estimates that a patient one standard deviation above average on BMI has an increased response expectation of 23.8 points. The bootstrap analysis (3.16) provided a standard error of 3.48 for this estimate. Bootstrapping, jackknifing and cross-validation require us to repeat the original estimation procedure for different data sets, which is easier to do if you know what the original procedure actually was.

Our thanks go to the discussants for their thoughtful remarks, and to the Editors for the formidable job of organizing this discussion.

Department of Statistics
Stanford University
Sequoia Hall
Stanford, California 94305-4065
USA
e-mail: brad@stat.stanford.edu